\documentclass[10pt, final]{amsproc}
\usepackage[cp850]{inputenc}
\usepackage{amssymb}
\usepackage[mathscr]{eucal}
\usepackage{amscd}
\usepackage{amsmath}
\usepackage{setspace}
\usepackage{tikz}
\def\Lip{\operatorname{Lip}}

\theoremstyle{plain}

\newtheorem{prop}{Proposition}
\newtheorem{lemma}{Lemma}
\newtheorem{cor}{Corollary}

\theoremstyle{remark}

\def\sep{\operatorname{sep}}
\title{The Kottman constant for $\alpha$-H\"older maps}

\author{Jes\'us Su\'arez de la Fuente}
\address{Escuela Polit\'ecnica, Avenida de la Universidad s/n, 10071 C\'aceres, Spain.}
\email{jesus@unex.es}

\begin{document}



\begin{abstract} We investigate the role of the Kottman constant of a Banach space $X$ in the extension of $\alpha$-H\"older continuous maps for every $\alpha\in (0,1]$.
\end{abstract}
\maketitle
\section{Introduction}
If $X$ is an infinite-dimensional Banach space then the Kottman constant \cite{Ko} of $X$ is defined as 
$$\kappa(X):=\sup_{(x_n)\in B(X)} \sep(x_n)$$
where for any sequence we define $$ \sep(x_n)=\inf_{n\neq m} \|x_n-x_m\|,$$
and $B(X)$ denotes the unit ball of $X$. Such a constant was introduced by Kottman \cite{Ko}. 
A well known result of Elton and Odell \cite{EO} asserts that $\kappa(X)>1$ for every infinite-dimensional Banach space $X$. Let us introduce a second parameter associated to a Banach space $X$. We define the constant $\lambda_1(X,c_0)$ as the infimum of all $\lambda>0$ such that for every subset $M$ of $X$ every Lipschitz map $f:M\to c_0$ admits an extension $F:X\to c_0$ with $\Lip(F)\leq \lambda \Lip(f)$. Kalton proved in \cite[Proposition 5.8]{Kal} the following unexpected result
\begin{prop} For every infinite dimensional Banach space $X$
$$\kappa(X)= \lambda_1 (X,c_0).$$
\end{prop}
 The aim of this note is to observe that the proof of \cite[Proposition 5.8]{Kal} contains the natural extension for $\alpha$-H\"older maps, see Proposition \ref{c0}. As far as we know, the first time that the Kottman constant was linked with the extension of $\alpha$-H\"older maps was in the proof of a result of Lancien and Randrianantoanina \cite[Theorem 2.2]{LR}. Although the Kottman constant is not mentioned, its role during the proof is quite evident. We present also a proof of their result where the Kottman constant appears explicitly, see Proposition \ref{c02}.
 
Recall that given metric spaces $(X,d)$ and $(Y,\rho)$ we say a map $f:X\to Y$ is $\alpha$-H\"older for $\alpha\in (0,1]$ if 
$$\Lip_{\alpha}(f)=\sup \left\lbrace \frac{\rho(f(x),f(y))}{d(x,y)^{\alpha}} : x\neq y \right\rbrace<\infty. $$
Therefore for each $\alpha\in(0,1]$ we may define the constant $\lambda_{\alpha}(X,Y)$ as the infimum of all $\lambda>0$ such that for every subset $M$ of $X$ every $\alpha$-H\"older map $f:M\to Y$ admits an extension $F:X\to Y$ with $\Lip_{\alpha}(F)\leq \lambda \Lip_{\alpha}(f)$. The paper is organized as follows. Section \ref{uno} contains the natural extension of \cite[Proposition 5.8]{Kal} to $\alpha$-H\"older maps that is our main result. It also contains the quantitative version of the result of Lancien and Randrianantoanina. 
Section \ref{dos} deals with $L_q$ instead of $c_0$ and gives a lower bound for $\lambda_{\alpha}(X,L_q)$.
\section{The estimate for $c_0$}\label{uno}
To prove our main result that is Proposition \ref{c0} we need the following easy Lemma. 
\begin{lemma}\label{lemma}
$$\lambda_{\alpha}(X,c_0)\geq \kappa(X)^{\alpha},\;\;\;\;\;\alpha\in (0,1].$$
\end{lemma}
\begin{proof}Let us denote by simplicity $\kappa=\kappa(X)$. Given $\varepsilon>0$ find $(x_n)\in B(X)$ such that
$$\kappa-\varepsilon\leq \|x_n-x_m\|,\;\;\;n\neq m.$$
Define the function $$f(x_n)=\left(\kappa-\varepsilon\right)^{\alpha}e_{n},\;\;\;n\in \mathbb{N}$$
 that is $\alpha$-H\"older with constant 1. Pick any extension $F$ of $f$ and denote $z=F(0)\in c_0$:
 $$\left\|F(0)-F(x_n) \right\| \leq \Lip_{\alpha}(F).$$
Hence we have for each coordinate $n\in \mathbb N$
$$\left|z(n)-\left(\kappa-\varepsilon\right)^{\alpha}\right| \leq \Lip_{\alpha}(F),$$
and taking limit as $n$ tends to infinity we have 
$$ (\kappa-\varepsilon)^{\alpha}\leq \Lip_{\alpha}(F).$$
Since the 
extension $F$ is arbitrary we have $$\left(\kappa-\varepsilon\right)^{\alpha}\leq \lambda_{\alpha}(X,c_0),$$
and letting  $\varepsilon\to 0$ we are done.
\end{proof}
The hard part of the proof of our main result is to prove the reverse inequality and this is due to Kalton. Thus we are ready to prove the natural extension of \cite[Proposition 5.8]{Kal}.
\begin{prop}\label{c0} For every infinite dimensional Banach space $X$
$$\lambda_{\alpha}(X,c_0)=\kappa(X)^{\alpha},\;\;\;\;\;\alpha\in (0,1].$$
\end{prop}
\begin{proof} Let us observe that
\begin{equation}\label{eq0}
\lambda_{\alpha}(X,c_0)=\lambda_1(X^{\alpha},c_0),
\end{equation}
where as usual $X^{\alpha}$ denotes the snowflake of $X$, this is, the metric space $X^{\alpha}=\left(X,\|\cdot\|^{\alpha}\right)$.
Therefore the result follows using the same proof of Kalton's result \cite[Proposition 5.8]{Kal}. Let us check it for the sake of completeness. Let us prove that the metric space $X^{\alpha}$ admits Lipschitz extension with constant $\kappa^{\alpha}$ that by Lemma \ref{lemma} will be enough. Suppose to the contrary that this last claim is false and let us reach a contradiction. 
Denoting by simplicity $\kappa=\kappa(X)$, \cite[Theorem 5.1(ii)]{Kal} shows that there exists $a \in  X^{\alpha}$,  $\varepsilon>0$ and a sequence $(x_n)\in X^{\alpha}$ such that 
\begin{equation}\label{eq4}
\kappa^{\alpha}\|x_k-a\|^{\alpha} + \varepsilon < \|x_j-x_k\|^{\alpha}, \;\;\;j<k.
\end{equation}
Since $\kappa^{\alpha}>1$ we may suppose $(x_n)$ is bounded. Hence replacing $x_n$ by $x_n-a$ and rescaling, we may find for some $\varepsilon'>0$ a sequence $(x_n')\in B(X^{\alpha})$ (and thus $(x_n')\in B(X)$) such that 
\begin{equation}\label{eq}
\kappa^{\alpha}\|x'_k\|^{\alpha} + \varepsilon' < \|x'_j-x'_k\|^{\alpha}, \;\;\;j<k.
\end{equation}
Observe that the expression above gives that $(x_n')$ does not converge to $0$ in $X^{\alpha}$, otherwise:
$$0<\varepsilon'\leq \kappa^{\alpha}\|x'_k\|^{\alpha} + \varepsilon' < \|x'_j-x'_k\|^{\alpha}\to 0,$$
and we would reach a contradiction. In particular, $(x_n')$ does not converge to $0$ in $X$. 
Now observe that given $\varepsilon_1>0$, using the fact that a bounded sequence of real numbers has a convergent subsequence, it is not hard to check that one may find $a,b$ and an infinite subset $\mathbb N_1$ of $\mathbb N$ such that
\begin{itemize}
\item[(i)] $$0<a\leq \|x'_n\|\leq b\leq 1,\;\;\;n\in \mathbb N_1.$$
\item[(ii)] $$1-\varepsilon_1\leq \frac{a}{b}.$$
\end{itemize}
Indeed, if $r\neq 0$ is the limit point then for the fixed $\varepsilon_1>0$ one may take $a=r-\theta$ and $b=r+\theta$ for $\theta>0$ small enough. We need a special choice of $\varepsilon_1>0$. To this end, consider the function  $f(t)=\kappa^{\alpha}(1-t)^{\alpha}$ whose limit as $t\to 0$ is $\kappa^{\alpha}$. Therefore, by definition of limit, given $\frac{\varepsilon'}{2}$, there is some $t_0>0$ such that $$f(t_0)\geq \kappa^{\alpha} - \frac{\varepsilon'}{2}.$$
Pick $a,b$ and $\mathbb N_1$ as above for the choice $\varepsilon_1=t_0$. Now, rescaling $1/b^{\alpha}$ the expression (\ref{eq}) for $j,k\in \mathbb N_1$, one has
\begin{eqnarray*}
\left \| \frac{x'_j}{b}-\frac{x'_k}{b} \right \|^{\alpha} &\geq & \kappa^{\alpha}\frac{\left \|x'_k\right \|^{\alpha}}{b^{\alpha}} + \frac{\varepsilon'}{b^{\alpha}}\\
&\geq & \kappa^{\alpha}\frac{\left \|x'_k\right \|^{\alpha}}{b^{\alpha}}+ \varepsilon'\\
&\geq & \kappa^{\alpha}\frac{a^{\alpha}}{b^{\alpha}}+ \varepsilon'\\
&\geq & \kappa^{\alpha}(1-\varepsilon_1)^{\alpha}+ \varepsilon'\\
&\geq & \kappa^{\alpha} +\frac{\varepsilon'}{2},
\end{eqnarray*}
where the last inequality follows by our choice of $\varepsilon_1>0$. This is, for $j,k\in \mathbb N_1$ with $j<k$ one has:
\begin{equation}\label{eq2}
\left \| \frac{x'_j}{b}-\frac{x'_k}{b} \right \|^{\alpha} \geq \kappa^{\alpha} +\frac{\varepsilon'}{2}.
\end{equation}
To finish observe that $$\left( \kappa^{\alpha} +\frac{\varepsilon'}{2} \right)^{\frac{1}{\alpha}}>\kappa$$ and pick $\rho>0$ such that 
$$\left (\kappa^{\alpha} +\frac{\varepsilon'}{2} \right)^{\frac{1}{\alpha}}-\kappa > \rho>0.$$
For this $\rho>0$, we have for $j,k\in \mathbb N_1$ taking $\alpha^{-1}$-power in (\ref{eq2}):
$$\left \| \frac{x'_j}{b}-\frac{x'_k}{b} \right \| \geq \kappa + \rho,\;\;\;\;j<k.$$
Since the points $(b^{-1}x'_n)_{n\in \mathbb N_1}\in B(X)$ we have reached a contradiction with the Kottman constant of $X$. 
Thus by \cite[Theorem 5.1]{Kal},  $X^{\alpha}$ has the Lipschitz $(\kappa^{\alpha},c_0)$-EP in the notation of \cite{Kal}. This is, $X^{\alpha}$ admits Lipschitz extension of $c_0$-valued Lipschitz maps with constant $\kappa^{\alpha}$. Hence $\lambda_1(X^{\alpha},c_0)\leq \kappa^{\alpha}$ and by Lemma \ref{lemma} and (\ref{eq0}) we are done.  Recall that  Kalton's argument shows also that the infimum defining $\lambda_{\alpha}(X,c_0)$ is attained.
\end{proof}
Recall that $\kappa(X)>1$, see \cite{EO}, and hence $\kappa(X)^{\alpha}>1$ for every $\alpha\in (0,1]$. In other words, Lemma \ref{lemma} shows there is no infinite-dimensional Banach space $X$ for which $\lambda_{\alpha}(X,c_0)=1$. This last was proved by Lancien and Randrianantoanina \cite[Theorem 2.2]{LR} replacing $c_0$ by a separable Banach space $Y$ containing an isomorphic copy of $c_0$. The next Proposition can be read as a quantitative version of \cite[Theorem 2.2]{LR} and the proof follows closely the original. Let us first introduce first some basic notation which can be found in \cite{K1}. We define a \textit{gauge} to be a function $\omega:[0,\infty)\to [0,\infty)$ which is a continuous increasing subadditive function satisfying $$\lim_{t\to 0}\omega(t)=0.$$ For the rest of this subsection $\omega$ will always denote a gauge. Recall that given metric spaces $(X,d)$ and $(Y,\rho)$, a map $f:X\to Y$ is $\omega$-Lipschitz if $$\Lip_{\omega}(f)=\sup \left\lbrace \frac{\rho(f(x),f(y))}{\omega(d(x,y))} : x\neq y \right\rbrace<\infty. $$ 
We may also define the constant $\lambda_{\omega}(X,Y)$ as the infimum of all $\lambda>0$ such that for every subset $M$ of $X$ every $\omega$-Lipschitz $f:M\to Y$ admits an extension $F:X\to Y$ with $\Lip_{\omega}(F)\leq \lambda \Lip_{\omega}(f)$.  Since $c_0$ is a $2$-absolute Lipschitz retract, it follows that $c_0$ is a 2-absolute $\omega$-Lipschitz retract. In other words, the extension of $\omega$-Lipschitz maps with values in $c_0$ is guaranteed with $\lambda_{\omega}(X,c_0)\leq 2$ for any Banach space $X$.  
\begin{prop}\label{c02} Let $Y$ be a separable Banach space containing an isomorphic copy of $c_0$. For every infinite dimensional Banach space $X$
$$\frac{\omega(\kappa(X))}{\omega(1)}\leq \lambda_{\omega}(X,Y).$$
\end{prop}
\begin{proof} Let us write for simplicity $\kappa=\kappa(X)$. Fix $\varepsilon >0$ and by James' distortion theorem \cite{J} pick an $(1+\varepsilon)$-isomorphic copy of $c_0$ in $Y$ through an into isomorphism $T$ with $\|T\|\leq 1 + \varepsilon$ and $\| T^{-1}\|\leq 1$, denote by $(e_n)$ the image of the canonical basis of $c_0$ by $T$ and by $(e_n^*)$ be the Hahn-Banach extensions to $Y$ of the corresponding coordinate functionals. By the separability assumption we may pick a subsequence $(e_{k_n}^*)$ of $(e_n^*)$ such that $(e_{k_n}^*)$ $w^*$-converges to some point $y^*\in Y^*$.  Given $\varepsilon_1>0$, pick $(x_n)\in B(X)$ for which $$\kappa-\varepsilon_1\leq \|x_n-x_m\|,\;\;\;\;\textit{for}\;\;n\neq m.$$ Define the map $f(x_n)=(-1)^{n}\omega\left(\kappa-\varepsilon_1\right)e_{k_n}$ for each $n\in \mathbb N$. We trivially find that $f$ is $\omega$-Lipschitz with constant less than or equal to $(1+\varepsilon)$.
Take any extension of $f$ to $X$, namely $F$, and observe that $$\|F(0)-F(x_n)\|\leq \Lip_{\omega}(F) \cdot\omega(\|x_n\|)\leq \Lip_{\omega}(F) \cdot\omega(1).$$ Therefore we find, for $y=F(0)$: $$|y(e_{k_n}^*)-(-1)^n\omega(\kappa-\varepsilon_1)|\leq \Lip_{\omega}(F)\cdot\omega(1),$$ so that $$-\Lip_{\omega}(F)\cdot\omega(1)-\omega(\kappa-\varepsilon_1)\leq y(e_{k_n}^*)\leq \Lip_{\omega}(F)\cdot\omega(1)-\omega(\kappa-\varepsilon_1),\;\;\:\;\textit{for n odd}.$$
$$-\Lip_{\omega}(F)\cdot\omega(1)+\omega(\kappa-\varepsilon_1)\leq y(e_{k_n}^*)\leq \Lip_{\omega}(F)\cdot\omega(1)+\omega(\kappa-\varepsilon_1),\;\;\;\:\textit{for n even}.$$ If $\Lip_{\omega}(F)\cdot\omega(1)<\omega(\kappa-\varepsilon_1)$, write $$\eta=\omega(\kappa-\varepsilon_1)-\Lip_{\omega}(F)\cdot\omega(1)>0.$$ From above we find that for $n$ even $y(e_{k_n}^*)\geq \eta$ while for $n$ odd one has $y(e_{k_n}^*)\leq -\eta$. Hence the sequence $(e_{k_n}^*)$ is not $w^*$-converging that is absurd. Therefore $\omega(\kappa-\varepsilon_1)\leq \Lip_{\omega}(F)\cdot\omega(1)$. Since this last must hold for every extension $F$ of $f$, we find that $\omega(\kappa-\varepsilon_1)\leq (1+\varepsilon)\lambda_{\omega}(X,Y)\cdot\omega(1)$ and letting $\varepsilon_1 \to 0$ and $\varepsilon\to 0$ we are done.
\end{proof} 
The result of Lancien and Randrianantoanina \cite[Theorem 2.2]{LR} for $\alpha$-H\"older maps may be recovered by taking the gauges $\omega(t)=t^{\alpha}$ with $\alpha \in (0,1]$.
Recall that $\ell_{\infty}$ is a 1-absolute $\omega$-Lipschitz retract, see \cite[Lemma 1.1.]{BL}.
Thus $\lambda_{\omega}(X,\ell_{\infty})=1$ for every infinite dimensional Banach space $X$ while $\omega(\kappa(X))>\omega(1)$ since $\omega$ is increasing. In particular, the separability assumption of Proposition \ref{c02} cannot be removed.
Let us give one last application of our results. Let us introduce  $\mathcal B_C(X,c_0)$ as the set of those $\alpha$ such that any $\alpha$-H\"older function $f$ from a subset of $X$ with values in $c_0$ and $\Lip_{\alpha}(f)=K$ can be extended to a function $F$ on the whole $X$ with $\Lip_{\alpha}(F)\leq C\cdot K$. Then Proposition \ref{c0} and a routine argument immediately gives
\begin{cor}For every infinite dimensional Banach space $X$
$$\mathcal B_C(X,c_0)=\left(0, \frac{  \log C} {\log \kappa (X)} \right],\;\;\;\;1\leq C\leq \kappa (X).$$
Or symmetrically, $$\mathcal B_{\kappa(X)^{\alpha}}(X,c_0)=\left(0, \alpha\right],\;\;\;\;\alpha\in (0,1].$$
\end{cor}
\section{An estimate for $L_q$}\label{dos}
Let us notice that the main idea of the proof of Lemma \ref{lemma} works replacing $c_0$ by many classic sequence spaces such as $\ell_q$-spaces. In general, the extension of $\alpha$-H\"older maps with values in $\ell_q$ is no guaranteed. However, for those cases in which there is extension, this is $\lambda_{\alpha}(X,\ell_q)<\infty$, the following bound could be useful
\begin{cor}\label{lq} For every infinite dimensional Banach space $X$
$$2^{-\frac{1}{q}}\kappa(X)^{\alpha}\leq \lambda_{\alpha}(X,\ell_q),\;\;\;\;\;\alpha\in (0,1].$$
\end{cor}
%
%
%
%
The situation for $L_q$ seems to be different so let us study a lower bound for $\lambda_{\alpha}(X,L_q)$. Since in many cases $\lambda_{\alpha}(X,L_q)=\infty$, it is clear that our lower bound only makes sense for those Banach spaces $X$ for which $\lambda_{\alpha}(X,L_q)<\infty$. Recall that Naor proves in \cite[Theorem 1]{N} that $\lambda_{\alpha}(L_p,L_q)<\infty$ for $\alpha\leq \frac{p}{2}$ and $1<p,q\leq 2$. Therefore our main motivation is the case of $L_q$ with $1<q\leq2$. The main technical obstruction to give a lower bound using the argument of Lemma \ref{lemma} is that there are no natural coordinates in $L_q$. To dodge this obstacle, we use a technical lemma due to Naor, see \cite[Lemma 2]{N}. For every $N\in \mathbb N$, put $\Omega=\left\lbrace  -1,+1 \right\rbrace ^N$ endowed with the uniform probability measure on $\Omega$ and denote by $r_1,...,r_N$ the Rademacher functions on $\Omega$.
\begin{lemma}(Naor)\label{naor}\\
For all $1<q<\infty$ and $Z\in L_q(\Omega)$
$$\frac{1}{N}\sum_{n=1}^N \mathbb E|Z-r_n|^q\geq 1-C\sqrt {\frac{\log N}{N}},$$ 
where $C$ depends only on $q$.
\end{lemma}
We are ready to present the main result of this section.
\begin{prop}\label{Lp} For every infinite-dimensional Banach space $X$ and $1<q \leq 2$
$$2^{-\frac{1}{q^*}}\kappa(X)^{\alpha}\leq \lambda_{\alpha}(X,L_q),\;\;\;\;\;\alpha\in (0,1].$$
\end{prop}
\begin{proof}
Let $\kappa,\lambda_{\alpha}$ be as in Proposition \ref{c0} and for given $\varepsilon>0$, pick $(x_n)\in B(X)$ for which $$\kappa-\varepsilon\leq \|x_n-x_m\|,\;\;\;\;\textit{for}\;\;n\neq m.$$ Fix $N\in \mathbb N$ and define the map $f_N(x_n)=\left(\kappa-\varepsilon\right)^{\alpha}r_n$ for each $n\in \left\lbrace 1,...,N \right\rbrace $. It turns out that $f_N$ is $\alpha$-H\"older with constant less than or equal to $2^{\frac{1}{q^*}}$
$$\|f_N(x_n)-f_N(x_m)\|_{L_q}=\left(\kappa-\varepsilon\right)^{\alpha}\|r_n-r_m\|_{L_q}\leq 2^{1-\frac{1}{q}}\|x_n-x_m\|^{\alpha}.$$
Take any extension of $f_N$ to $X$, namely $F$, and observe that
$$\|F(0)-F(x_n)\|_{L_q}\leq \Lip_{\alpha}(F) \|r_n\|^{\alpha}\leq \Lip_{\alpha}(F).$$
Let us denote $F(0)=Z\in L_q(\Omega)$,  we have that 
\begin{eqnarray*}
\Lip_{\alpha}(F)^q&\geq & \frac{1}{N}\sum_{n=1}^N\mathbb E|Z-(\kappa-\varepsilon)^{\alpha}r_n|^q\\
&=& (\kappa-\varepsilon)^{\alpha q}\frac{1}{N}\sum_{n=1}^N\mathbb E|(\kappa-\varepsilon)^{-\alpha}Z-r_n|^q\\
&\geq&  (\kappa-\varepsilon)^{\alpha q}\left(1-C\sqrt {\frac{\log N}{N}}\right),
\end{eqnarray*}
where the last inequality follows from Lemma \ref{naor}. Since this must hold for every extension $F$ of $f_N$ and every $N\in \mathbb N$ we find that
$$(\kappa-\varepsilon)^{\alpha}\leq \lambda_{\alpha}2^{\frac{1}{q^*}},$$
and letting  $\varepsilon\to 0$ we are done.
\end{proof}
Since it is well known that $\kappa(L_p)=2^{\frac{1}{p}}$ for $1<p\leq 2$, Proposition \ref{Lp} yields
\begin{cor}\label{c} For $1<p,q\leq 2$,
$$2^{\frac{\alpha}{p}-\frac{1}{q^*}}\leq \lambda_{\alpha}(L_p,L_q).$$
\end{cor}
Recall that Naor proves also in \cite[Theorem 1]{N} that there is no isometric extension for $\alpha>\frac{p}{q^*}$ and $1<p,q\leq 2$. 
As Corollary \ref{c}  shows, the Kottman constant explains geometrically why there is no isometric extension for these values. To finish, Corollary \ref{c} gives $$1\leq \lambda_{\alpha}(L_p,L_q)$$ for $\alpha=\frac{p}{q^*}$ and $1<p,q\leq 2$ while \cite[Theorem 1]{N} provides us with $\lambda_{\alpha}(L_p,L_q)=1$, so it is sharp. 

\end{document}